\documentclass[11pt,leqno]{amsart}
\usepackage{amssymb,verbatim,enumerate,ifthen,hyperref}
\usepackage{mathtools} 

\usepackage{tikz-network}
\usepackage{graphicx}
\usepackage{mathabx}
\usepackage{accents}
\usepackage[mathscr]{eucal}
\usepackage[utf8]{inputenc}
\usepackage[T1]{fontenc}
\def\R{\mathbb{R}}

\def\H{\mathscr{H}}

\long\def\comment#1{}

\newtheorem{theorem}{Theorem}[section]
\newtheorem*{theorem*}{Theorem}
\def\Thm#1#2{\ifthenelse{\equal{#1}{*}}{\begin{theorem*}#2\end{theorem*}}
             {\begin{theorem}\label{T#1}#2\end{theorem}}}
\newtheorem{Atheorem}{Theorem}

\newtheorem{proposition}[theorem]{Proposition}
\newtheorem*{proposition*}{Proposition}
\def\Prp#1#2{\ifthenelse{\equal{#1}{*}}{\begin{proposition*}#2\end{proposition*}}
{\begin{proposition}\label{P#1}#2\end{proposition}}}
\def\prp#1{Proposition~\ref{P#1}}

\newtheorem{corollary}[theorem]{Corollary}
\newtheorem*{corollary*}{Corollary}
\def\Cor#1#2{\ifthenelse{\equal{#1}{*}}{\begin{corollary*}#2\end{corollary*}}
             {\begin{corollary}\label{C#1}#2\end{corollary}}}

\newtheorem{lemma}[theorem]{Lemma}
\newtheorem*{lemma*}{Lemma}
\def\Lem#1#2{\ifthenelse{\equal{#1}{*}}{\begin{lemma*}#2\end{lemma*}}
             {\begin{lemma}\label{L#1}#2\end{lemma}}}

\theoremstyle{definition}
\newtheorem{remark}[theorem]{Remark}
\newtheorem*{remark*}{Remark}
\def\Rem#1#2{\ifthenelse{\equal{#1}{*}}{\begin{remark}\rm #2\end{remark}}
             {\begin{remark}\label{R#1}\rm #2\end{remark}}}

\newtheorem{example}[theorem]{Example}
\newtheorem*{example*}{Example}
\def\Exa#1#2{\ifthenelse{\equal{#1}{*}}{\begin{example*}\rm #2\end{example*}}
             {\begin{example}\label{Ex#1}\rm #2\end{example}}}

\def\eq#1{{\rm(\ref{E#1})}}
\def\Eq#1#2{\ifthenelse{\equal{#1}{*}}
  {\begin{equation*}\begin{aligned}#2\end{aligned}\end{equation*}}
  {\begin{equation}\begin{aligned}\label{E#1}#2\end{aligned}\end{equation}}}

\begin{document}
\begin{flushright}
\end{flushright}
\vspace{5mm}

\date{\today}

\title[On Weighted Monotone and Subadditive Graphs]
{On Weighted Monotone and Subadditive Graphs}

\author[A. R. Goswami]{Angshuman R. Goswami}
\address[A. R. Goswami]{Department of Mathematics, University of Pannonia,
H-8200 Veszprem, Hungary}
\email{\{goswami.angshuman.robin@mik.uni-pannon.hu}


\thanks{The research of the author was supported by the `Ipar a Veszprémi Mérnökképzésért' foundation.
}

\begin{abstract}
Let $G(V,E)$ be a graph, and $\H:=\big\{H:H\subseteq G\big\}$ denote the collection of all possible subgraphs of $G$. Then for each non-negative function $w:\H\to\R_+$, the graph $G(V,E,w)$ is said to be a weighted graph.
\\

A weighted graph $G(V,E,w)$ is called monotone (increasing), if for any $H_1,H_2\subseteq G$ with 
$H_1\subset H_2$, the following inequality holds:
\Eq{*}{
w\big(H_1\big)\leq w\big(H_2\big).
}

On the other hand, a weighted graph $G(V,E,{w})$ is termed subadditive, if for any 
$H_1,H_2\subseteq G$, the following discrete functional inequality is satisfied:
\Eq{*}{
{w}\big(H_1\cup H_2\big)\leq {w}\big(H_1\big)+ {w}\big(H_2\big).
}

Our main result demonstrates that for any graph $G(V,E,w)$, it is possible to construct both the largest monotone and the greatest subadditive minorants. In other words, it is feasible to formulate the largest increasing function $\overline{w}:\H\to\R_+$ and subadditive function $\widetilde{w}:\H\to\R_+$ such that $\overline{w}(H)\leq w(H)$ and $\widetilde{w}(H)\leq w(H)$ hold respectively for all $H\subseteq G$ .
\end{abstract}

\maketitle
\section*{Introduction}
Throughout this paper, the $\R_+$ will denote the set of non-negative numbers. The symbol $G(V,E)$ is used to refer to the non-weighted graph with vertex and edge sets $V$ and $E$, respectively. Let $\H:=\big\{H:H\subseteq G\big\}$ be the class that consists of all possible subgraphs of $G$. The notation $G(V,E,w)$ represents the weighted version of the graph $G(V,E)$ after assigning the weight $w:\H\to\R_+$ to it. Additionally, it is worth noting that the graphs $G(V,E,w_1)$ and $G(V,E,w_2)$ share the same vertex and edge sets, but their weights differ.
\\

The notions of monotonicity and subadditivity in the context of graph theory are not new. Mathematicians have introduced several definitions of monotone and subadditive graphs for use in various types of discrete analysis.\\

A weighted a graph $G(V,E,w)$ is said to be \textit{monotone (increasing)}, if the following inequality holds along with mentioned conditions:
\Eq{*}{
w\big(H_1\big)\leq w\big(H_2\big)\quad \mbox{for any}\quad H_1,H_2\subseteq G \quad\mbox{with}\quad H_1\subset H_2.
}
The chromatic number $\chi(G)$, the clique number $\omega(G)$, and the maximum degree $\Delta(G)$ are some of the parameters that are closely aligned to our definition of monotonicity. Some works in this direction can be found in 
\cite{Marko, Saleh, Dirac} and references therein.\\

Similarly, a weighted graph $G(V,E,w)$ is said to be \textit{subadditive}, if for any $H_1,H_2\subseteq G$, the following discrete functional inequality satisfied:
\Eq{*}{w\big(H_1\cup H_2\big)\leq w\big(H_1\big)+w\big(H_2\big).} For more details, one can look into \cite{Rezaei, Zhao, Schaefer} and references.\\

Let $f : X \to \R$ be a real-valued function. 
A function $g : X \to \R$ that satisfies a specific property $P$ (such as monotonicity, convexity, subadditivity, periodicity, continuity, differentiability, etc.) is called a \textit{$P$-{minorant}} of $f$ if 
$$
g(x) \leq f(x) \quad \text{holds for all} \quad x \in X.
$$
In case of reverse inequality, $g$ is labelled as \textit{$P$-majorant}.
Among all $P$-minorants, we may look for the largest one — called the \textit{greatest(largest) $P$-minorant} of $f$. Similarly, one can also define the smallest \textit{$P$-majorant}.\\

In the classical function theory, for any function $f:I(\subseteq\R)\to\R$, we have explicit formulas to derive various $P$-minorants. For instance, in Theorem 3 of \cite{Zsolt} and in Proposition 3.3 of \cite{Goswami}, we can observe the constructions of the largest increasing and subadditive minorants, respectively. The discrete analogues of these functional concepts can be found in recent papers \cite{Goswamii, Goswamiii}. Moreover, it is worth mentioning that some related ideas of monotone minorants, though in a different context, appear in the papers \cite{Gord, Sinnamon, Anantharam}. All these offer significant insights into the topic, including the poset-based monotone envelopes. \\

These findings motivate us to formulate the largest $P$-minorant and the smallest $P$-majorant in the framework of graphs as well. For any given graph $G(V,E,w)$, we show how to formulate the largest/smallest monotone minorants. We also propose a formula to obtain the greatest subadditive minorant. We present similar results for the case of majorants as well. Besides, we study some interesting overlaps of these two properties.\\

We start our investigations with the monotone minorants.
\section{Main Results}
Throughout the section, we assume that the graph $G(V,E)$ has a finite number of vertices and edges. The results are presented accordingly. However, one can easily transform the findings to infinite settings.

\Prp{1} {Let $G(V,E,w)$ be a weighted graph. The function $\overline{w}:\H\to\R_+$ is defined as follows:
\Eq{1}{
\overline{w}\big(H\big):=\min\Big\{w\big(H'\big)\,\,\Big|\,\,H\subseteq H'\Big\}.
}
Then $G(V,E,\overline{w})$ is the largest monotonically increasing minorant of $G(V,E,w)$.
}
\begin{proof}
By definition, we have 
\Eq{1}{
\overline{w}\big(H\big)=\min\Big\{w\big(H'\big)\,\,\Big|\,\,H\subseteq H'\Big\}\leq w(H).
}
This shows that $G(V,E,\overline{w})$ is a minorant of $G(V,E,{w})$.\\

To prove the monotonicity, we assume $H_1,H_2\subseteq G$ such that $H_1\subset H_2$ holds. Then using \eq{1}, we can compute the following
\Eq{*}{
\overline{w}\big(H_1\big)= \min\Big\{w\big(H'\big)\,\,\Big|\,\,H_1\subseteq H'\Big\} \leq \min\Big\{w\big(H'\big)\,\,\Big|\,\,H_2\subseteq H'\Big\}=\overline{w}\big(H_2\big).
}
This establishes the increasing monotonicity of $G(V,E,\overline{w})$.\\

If possible, let there exist another mapping $w_0:\H\to\R_+$ such that $G(V,E,w_0)$ is the largest increasing minorant of $G(V,E,w)$. This together with the monotonicity of $G(V,E,w_0)$ implies that for any $H\subset G$, all the existing $H'\subseteq G$ with $H\subseteq H'$, satisfy the following inequality
\Eq{*}{
\overline{w}\big(H\big)\leq{w_0}\big(H\big)\leq w_0\big(H'\big)\leq w\big(H'\big).
}
Now, taking the minimum to the right-most side of the above inequality, we have 
\Eq{*}{
{w_0}\big(H\big)\leq \min\Big\{w\big(H'\big)\,\,\Big|\,\,\mbox{for all}\,\,H\subseteq H'\Big\}=\overline{w}\big(H\big).
}
This contradicts our assumption on $G(V,E,w_0)$.\\

Hence, $G(V,E,\overline{w})$ is the largest increasing minorant of $G(V,E,w)$. This completes the proof.
\end{proof}
Based on the above result, the following corollary can be stated. The establishment of it is straightforward. Hence, we only provide the sketch of the proof.
\Cor{1}{Let $G(V,E,w_1)$ and $G(V,E,w_2)$ are two weighted graphs such that for any 
$H\subset G$ the inequality $w_1\big(H\big)\leq w_2\big(H'\big)$ holds for all $H'\subseteq G$ with $H\subseteq H'.$ Then there exists a monotone (increasing) graph $G(V,E,\overline{w})$, that satisfies the following discrete functional inequality
\Eq{*}{
w_1\big(H\big)\leq \overline{w}\big(H\big)\leq w_2\big(H\big)\qquad\mbox{for all}\qquad H\subseteq G.
}
}
\begin{proof}
We construct the weight function $\overline{w}:\H\to\R_+$ as in \eq{1}. From construction, \newline $w_1\big(H\big)\leq \overline{w}\big(H\big)\leq w_2\big(H\big)$ is evident. The increasing monotonicity of $G(V,E,\overline{w})$ can be shown as in \prp{1}.
\end{proof}
The proposition below is analogous to \prp{1}. Hence, the proof is left for the readers.
\Prp{11} {Let $G(V,E,w)$ be a weighted graph. The function $\overline{w}:\H\to\R_+$ is defined as follows:
\Eq{*}{
\underline{w}\big(H\big):=\max\Big\{w\big(H'\big)\,\,\Big|\,\,H'\subseteq H\Big\}.
}
Then $G(V,E,\underline{w})$ is the smallest monotone(increasing) majorant of $G(V,E,w)$.
}
The following theorem provides the methodology to obtain the largest subadditive minorant for any given graph. 
\Prp{2}{ Let $G(V,E,w)$ be a weighted graph. The function $\widetilde{w}:\H\to\R_+$ is defined as follows:
\Eq{3}{
\widetilde{w}\big(H\big):=\min\Big\{w\big(H_1\big)+\cdots +w\big(H_n\big)\,\,\Big|\,\,H_1,\cdots,H_n\subseteq H\,\,\mbox{such that}\,\,\overset{n}{\underset{i=1}{\cup}}H_i=H \}.
}
Then $G(V,E,\widetilde{w})$ be the largest subadditive minorant of $G(V,E,{w})$.
}
\begin{proof}
From the definition, it is clearly visible that 
\Eq{*}{
\widetilde{w}\big(H\big)=\min\Big\{w\big(H_1\big)+\cdots +w\big(H_n\big)\,\,\Big|\,\,H_1,\cdots,H_n\subseteq H\,\,\mbox{such that}\,\,\overset{n}{\underset{i=1}{\cup}}H_i=H \}\leq w\big(H\big).}
This shows that $G(V,E,\widetilde{w})$ is a minorant of $G(V,E,{w})$.\\

To prove subadditivity, we assume two distinct arbitrary subgraphs $H_1$ and $H_2$ of $G$. Then from \eq{3}, we have the following
\Eq{4}{
\widetilde{w}(H_1)=w\big(H_1^{1}\big)+\cdots +w\big(H_1^{m}\big)\,\,&\mbox{such that}\,\,H_1^{1},\cdots,H_1^{m}\subseteq H_1\,\,\mbox{satisfying}\,\,\overset{m}{\underset{i=1}{\cup}}H_1^{i}=H_1\\
&\quad\,\,\mbox{and}\\
\widetilde{w}(H_2)=w\big(H_2^{1}\big)+\cdots +w\big(H_2^{n}\big)\,\,&\mbox{such that}\,\,H_2^{1},\cdots,H_2^{n}\subseteq H_2\,\,\mbox{satisfying}\,\,\overset{n}{\underset{i=1}{\cup}}H_2^{i}=H_2.
}
From the above, we can also observe the following inclusion
\Eq{*}{H_1^1\cup\cdots\cup H_1^m\cup H_2^1\cup\cdots\cup H_2^n=H_1\cup H_2\subseteq G.
}
Thus from \eq{3}, we can conclude the following
\Eq{*}{
w\big(H_1\cup H_2\big)\leq w\big(H_1^{1}\big)+\cdots +w\big(H_1^{m}\big)+w\big(H_2^{1}\big)+\cdots +w\big(H_2^{n}\big).
}
Using \eq{4} in the above inequality, we finally have 
\Eq{*}{w\big(H_1\cup H_2\big)\leq w\big(H_1\big)+w\big(H_2\big).}
Since, $H_1,H_2\subseteq G$ are arbitrary, hence this yields that the weighted graph $G(V,E,\widetilde{w})$ possesses subadditivity.\\

Now, if possible, let there exist another mapping $w_0:\H\to\R_+$ such that $G(V,E,w_0)$ is the largest subadditive minorant of $G(V,E,w)$. This together with the subadditivity of $G(V,E,w_0)$ implies the following inequality
\Eq{*}{
w_0\big(H\big)&\leq w_0\big(H_1\big)+\cdots +w_0\big(H_n\big)\,\,\mbox{for all}\,\,H_1,\cdots,H_n\subseteq H\,\,\mbox{satisfying}\,\,\overset{n}{\underset{i=1}{\cup}}H_i=H\\
&\leq w\big(H_1\big)+\cdots +w\big(H_n\big).
}
Taking the minimum to the right-most part of the above inequality, we have
\Eq{*}{
w_0\big(H\big)
\leq\min\Big\{ w\big(H_1\big)+\cdots +w\big(H_n\big)\,\,\Big|\,\,H_1,\cdots,H_n\subseteq H\,\,\mbox{provided}\,\,\overset{n}{\underset{i=1}{\cup}}H_i=H\Big\}=\widetilde{w}\big(H\big).
}
This contradicts our assumption on $G(V,E,w_0)$.\\

Hence, $G(V,E,\widetilde{w})$ is the largest subadditive minorant of $G(V,E,w)$. This completes the proof.
\end{proof}
In the next proposition, we will observe the close interaction between graphical monotone and subadditive properties.
\Prp{3}{In the above theorem, if $G(V,E,w)$ is monotone then $G(V,E,\widetilde{w})$ also possesses same monotonicity.}
\begin{proof}
To prove the assertion, first assume that $G(V,E,w)$ is monotonically increasing. Let $H_1,H_2\subseteq G$ such that $H_1\subset H_2$. Then from \eq{3} of \prp{2}, there must exists subgraphs $H_1^{1},\cdots,H_1^{m}\subseteq H_1$ and $H_2^{1},\cdots,H_2^{n}\subseteq H_2$ such that both equalities of \eq{4} holds along with the mentioned conditions. Also, one can easily observe the following inclusions
\Eq{*}{
H_2^{1}\cap H_1\subset H_1,\cdots, H_2^{n}\cap H_1\subset H_1 \quad\mbox{which also implies}\quad \overset{n}{\underset{i=1}{\cup}}\big(H_2^{i}\cap H_1\big)=H_1.
}
Using these and then utilizing the monotonicity of $G$, the first equality of \eq{4} can be extended as follows
\Eq{*}{
\widetilde{w}(H_1)&=w\big(H_1^{1}\big)+\cdots +w\big(H_1^{m}\big)\\
&\leq w\big(H_2^{1}\cap H_1\big)+\cdots +w\big(H_2^{n}\cap H_1\big)\\
&\leq  w\big(H_2^{1}\big)+\cdots +w\big(H_2^{n}\big)\\
&=\widetilde{w}(H_2).
}
Since, $H_1,H_2$ are arbitrary, the above inequality establishes increasing monotonicity of $G(V,E,\widetilde{w})$.\\

Now, we assume $G(V,E,w)$ is monotonically decreasing and $H_1,H_2$ are two arbitrary subsets of it such that $H_1\subset H_2$ holds. The decreasing monotonicity of $G$ implies the following inequalities
\Eq{*}{
w\big(H_1\cup H_2\big)\leq w\big(H_1\big)\qquad \mbox{and}\qquad w\big(H_1\cup H_2\big)\leq w\big(H_2\big).
}
This implies, $w\big(H_1\cup H_2\big)\leq w\big(H_1\big)+w\big(H_2\big)$, which is obvious. That is, $G(V,E,{w})$ possesses subadditivity. In other words,
if $G(V,E,{w})$ is decreasing, then $G(V,E,{w})=G(V,E,\widetilde{w})$.\\

This validates the assertion and completes the proof.
\end{proof}
The following corollary is a direct consequence of \prp{2}. Therefore, we provide only the draft of the proof.
\Cor{2}{Let $G(V,E,w_1)$ and $G(V,E,w_2)$ are two weighted graphs such that 
\Eq{*}{w_1\big(H\big)\leq w_2\big(H_1\big)+\cdots+w_2\big(H_n\big)
}
holds for any $H_1,\cdots,H_n\subseteq H$ with $\overset{n}{\underset{i=1}{\cup}}H_i=H$.
Then there exists a subadditive graph $G(V,E,\widetilde{w})$, that satisfies the following discrete functional inequality
\Eq{*}{
w_1\big(H\big)\leq \widetilde{w}\big(H\big)\leq w_2\big(H\big)\qquad\mbox{for all}\qquad H\subseteq G.
}
}
\begin{proof}
We construct the weight function $\widetilde{w}:\H\to\R_+$ as in \eq{3}. From there, the inequality $w_1\big(H\big)\leq \widetilde{w}\big(H\big)\leq w_2\big(H\big)$, for all $H\subseteq G$ is obvious. The proof of the remaining assertions related to subadditivity are analogous to \prp{2}.
\end{proof}
The next result is similar to \prp{2}. Thus, we decided to leave the proof to the readers.
\Prp{22}{ Let $G(V,E,w)$ be a weighted graph. The function $\underset{\widetilde{\,\,\,\,\,\,\,}}{w}:\H\to\R_+$ is defined as follows:
\Eq{*}{
\underset{\widetilde{\,\,\,\,\,\,\,}}{w}\big(H\big):=\max\Big\{w\big(H_1\big)+\cdots +w\big(H_n\big)\,\,\Big|\,\,H_1,\cdots,H_n\subseteq H\,\,\mbox{such that}\,\,\overset{n}{\underset{i=1}{\cup}}H_i=H \}.
}
Then $G(V,E,\underset{\widetilde{\,\,\,\,\,\,\,}}{w})$ be the smallest super-additive majorant of $G(V,E,{w})$.
}
Most of the results presented in this paper can be naturally extended to infinite graphs by replacing $\min$ and $\max$ with $\inf$ and $\sup$, respectively at the appropriate points. Furthermore, the paper introduces several computational aspects related to these theoretical findings. In particular, it is possible to compute the largest/smallest monotones and the largest subadditive minorants for a given weighted graph $G(V,E,w)$ when the number of vertices and edges remains within a manageable range. However, for graphs of substantial size, the computation of these specific minorants may becomes computationally intractable. To tackle this challenge, future research may focus on developing algorithms capable of approximating these minorants for large-scale graphs.
\bibliographystyle{plain}

\end{document}